\documentclass[10pt]{article} \usepackage{amsfonts} \usepackage{amssymb}\usepackage{amsmath}\usepackage{amsfonts}
\usepackage[T2A]{fontenc}
\usepackage[cp1251]{inputenc}
\usepackage[english,russian]{babel}\usepackage[T2A]{fontenc}
\usepackage[cp1251]{inputenc}
\usepackage[english,russian]{babel}

\makeatletter
\renewcommand{\@makefnmark}{}
\makeatother
\newcommand{\elx}{e^{i\lambda x}}
\newcommand{\elxx}{e^{-i\lambda x}}

\newcommand{\elxb}{e^{i\bar\lambda x}}
\newcommand{\elxxb}{e^{-i\bar\lambda x}}

\begin{document}
\baselineskip=10pt
\pagestyle{plain}
{\Large

\medskip
\medskip
\medskip

\footnote{

\hspace{2mm}Keywords: Dirac system, maximal operator, spectral expansions}

\centerline {\bf On a necessary condition of convergence of spectral expansions}
 \centerline {\bf corresponding to the Dirac operators}

\medskip
\medskip
\medskip
\medskip
\medskip

\centerline { Alexander Makin}
\medskip
\medskip
\medskip

\medskip
\medskip
\begin{quote}{\normalsize
The paper is concerned with the basis properties of root function systems of  the Dirac operator with a complex-valued summable potential. We establish  a necessary condition of convergence of  corresponding spectral expansions.}
   \end{quote}

\centerline {\bf 1. Introduction}

\medskip
\medskip
\medskip
\medskip

In the  present paper, we study the Dirac system
$$
B\mathbf{y}'+V\mathbf{y} =\lambda\mathbf{y},\eqno(1)
$$
 where $\mathbf{y}={\rm col}(y_1(x),y_2(x))$,
 \[
 B=\begin{pmatrix}
 -i&0\\
 0&i
 \end{pmatrix},\quad V=\begin{pmatrix}
 0&P(x)\\
 Q(x)&0
 \end{pmatrix},
\]
 and
the functions $P, Q\in L_1(0,\pi)$
in the space
$\mathbb{H}=L_2(0,\pi)\oplus L_2(0,\pi)$ with the inner product

$$
\langle\mathbf{f},\mathbf{g}\rangle=\int_0^\pi(f_1(x)\overline{g_1(x)}+f_2(x)\overline{g_2(x)})dx,
$$
where $\mathbf{f}={\rm col}(f_1(x),f_2(x)),\quad\mathbf{g}={\rm col}(g_1(x),g_2(x))$.

Let us define in the space
$\mathbb{H}$ the maximal operator
$\mathbb{L}\mathbf{y}=B\mathbf{y}'+V\mathbf{y}$
with the domain $D(\mathbb{L})=\{\mathbf{y}\in W_1^1[0,\pi]\oplus  W_1^1[0,\pi]:\, \mathbb{L}\mathbf{y}\in \mathbb{H}$.

By an eigenfunction of operator $\mathbb{L}$ corresponding to a complex eigenvalue $\lambda$ we mean any complex-valued function $\stackrel{0}{\mathbf{y}}\in D(\mathbb{L})$ $( \stackrel{0}{\mathbf{y}}\not\equiv0$) satisfying the equation $\mathbb{L}\stackrel{0}{\mathbf{y}}=\lambda\stackrel{0}{\mathbf{y}}$ almost everywhere in the domain $G=(0,\pi)\times(0,\pi)$. Analogously, by an associated function $\stackrel{l}{\mathbf{y}}$ of operator $\mathbb{L}$ of order $l=1,2,\ldots$ corresponding to the same $\lambda$
and the  eigenfunction $\stackrel{0}{\mathbf{y}}$  mean any complex-valued function $\stackrel{l}{\mathbf{y}}\in D(\mathbb{L})$
satisfying the equation $\mathbb{L}\stackrel{l}{\mathbf{y}}=\lambda\stackrel{l}{\mathbf{y}}+\stackrel{l-1}{\mathbf{y}}$ almost everywhere in  $G$.
Thus, the root functions are understood in the generalized sense (irrespective of the boundary conditions).

Let $\{\mathbf{y}_n\}$ $(n\in\mathbb{N})$ be a root function system of operator $\mathbb{L}$. The main aim of the present article is to find necessary conditions for the basis property
in $\mathbb{H}$ for the indicated system.

\medskip
\medskip
\medskip
\centerline {\bf 2. Main results}
\medskip
\medskip
\medskip

Denote by $\mathbb{L}^*$ the operator formally adjoint of the operator $\mathbb{L}$:

$$
\mathbb{L}^*=B\mathbf{y}'+V^*\mathbf{y},
$$
where
\[
V^*=\begin{pmatrix}
 0&\bar Q(x)\\
 \bar P(x)&0
 \end{pmatrix}.
\]

Together with system (1) we will consider the system

$$
B\mathbf{z}'+V^*\mathbf{z} =\bar\lambda\mathbf{z},\eqno(1^*)
$$
where $\mathbf{z}={\rm col}(z_1(x),z_2(x))$.

{\bf Lemma 1.} {\it Let $\mathbf{y}$ and $\mathbf{z}$ be solutions to equations (1) and (1*) correspondingly.
If
 $$
 |{\rm Im}\lambda|\to\infty
$$
 and
$$
\langle\mathbf{y},\mathbf{z} \rangle=1,
$$
then
$$
\|\mathbf{y}\|_\mathbb{H}\|\mathbf{z}\|_\mathbb{H}\to\infty.
$$
}

Proof. Denote by
$\Omega_+=\{{\rm Im}\lambda>-r, |\lambda|>\lambda_0\}$, $\Omega_-=\{{\rm Im}\lambda<r, |\lambda|>\lambda_0\}$,
$(r>0)$.
It is known [1, Th. 6.8] (see also [2-6]) that in the domain $\Omega_+$ $(\Omega_-)$ system (1) has a fundamental solution matrix

 \[
Y_\pm(x,\lambda)=\|Y_{\pm,jk}(x,\lambda)\|=\begin{pmatrix}
(1+b_{\pm,11}(x,\lambda))\elx& b_{\pm,12}(x,\lambda)\elxx \\
b_{\pm,21}(x,\lambda)\elx&(1+ b_{\pm,22}(x,\lambda))\elxx
\end{pmatrix},\eqno(2)
\]

where
$$
b_{\pm,jk}(x,\lambda)=o(1)\eqno(3)
$$
as $|\lambda|\to\infty$
in the domain $\Omega_+$ $(\Omega_-)$ uniformly in $x\in[0,\pi]$ and
$$
\|b_{\pm,jk}(x,\lambda)\|_{W_1^1[0,\pi]}<c_0.
$$

 Analogously, in the domain $\Omega_+$ $(\Omega_-)$ system (1*) has a fundamental solution matrix

 \[
Z_\mp(x,\lambda)=\|Z_{\mp,jk}(x,\bar\lambda)\|=\begin{pmatrix}
(1+\tilde b_{\mp,11}(x,\bar\lambda))\elxb& \tilde b_{\mp,12}(x,\bar\lambda)\elxxb \\
\tilde b_{\mp,21}(x,\bar\lambda)\elxb&(1+ \tilde b_{\mp,22}(x,\bar\lambda))\elxxb
\end{pmatrix},\eqno(2^*)
\]

where
$$
\tilde b_{\mp,jk}(x,\bar\lambda)=o(1)\eqno(3^*)
$$
as $|\lambda|\to\infty$
in the domain  $\Omega_+$ $(\Omega_-)$  uniformly in $x\in[0,\pi]$ and
$$
\|\tilde b_{\mp,jk}(x,\bar\lambda)\|_{W_1^1[0,\pi]}<c_0.
$$

Let
$Y^{[k]}(x,\lambda)$
be the $k$th column of matrix (2) and
$Z^{[k]}(x,\lambda)$
be the $k$th column of matrix (2*).
Evidently,
$$
\mathbf{y}=C_1Y^{[1]}+C_2Y^{[2]},\quad\mathbf{z}=\tilde C_1Z^{[1]}+\tilde C_2Z^{[2]},
$$
where $C_1,C_2,\tilde C_1,\tilde C_2$ are some complex numbers.
Consider the case $\lambda\in\Omega_+$.

It is easy to see that
$$
\frac{c_1}{\sqrt{|{\rm Im}\lambda}|+1}\le||e^{i\lambda t}||_{L_2(0,\pi)}\le\frac{c_{2}}{\sqrt{|{\rm Im}\lambda|}+1},\eqno(4)
$$

$$
\frac{c_1e^{\pi|{\rm Im}\lambda|}}{\sqrt{|{\rm Im}\lambda|}+1}\le||e^{-i\lambda t}||_{L_2(0,\pi)}\le\frac{c_{2}e^{\pi|{\rm Im}\lambda|}}{\sqrt{|{\rm Im}\lambda|}+1}.\eqno(5)
$$
Combining (2), (2*), (3), (3*), (4), (5), we obtain

$$
\frac{c_9}{\sqrt{|{\rm Im}\lambda}|+1}\le||Y_{+,11}(x,\lambda)||_{L_2(0,\pi)}\le\frac{c_{10}}{\sqrt{|{\rm Im}\lambda|}+1},\eqno(6)
$$

$$
\frac{c_9e^{\pi|{\rm Im}\lambda|}}{\sqrt{|{\rm Im}\lambda|}+1}\le||Y_{+,22}(x,\lambda)||_{L_2(0,\pi)}\le\frac{c_{10}e^{\pi|{\rm Im}\lambda|}}{\sqrt{|{\rm Im}\lambda|}+1},\eqno(7)
$$

$$
||Y_{+,21}(x,\lambda)||_{L_2(0,\pi)}=\frac{o(1)}{\sqrt{|{\rm Im}\lambda|}+1},\eqno(8)
$$

$$
||Y_{+,12}(x,\lambda)||_{L_2(0,\pi)}=\frac{o(1)e^{\pi|{\rm Im}\lambda|}}{\sqrt{|{\rm Im}\lambda|}+1}.\eqno(9)
$$

It follows from (6-9) that

$$
\frac{c_9}{\sqrt{|{\rm Im}\lambda}|+1}\le||Y^{[1]}||_\mathbb{H}\le\frac{c_{10}}{\sqrt{|{\rm Im}\lambda|}+1},\eqno(10)
$$

$$
\frac{c_9e^{\pi|{\rm Im}\lambda|}}{\sqrt{|{\rm Im}\lambda|}+1}\le||Y^{[2]}||_\mathbb{H}\le\frac{c_{10}e^{\pi|{\rm Im}\lambda|}}{\sqrt{|{\rm Im}\lambda|}+1}.\eqno(11)
$$

Inequalities  (10), (11)   imply
$$
|C_1|^2\|Y^{[1]}\|^2_\mathbb{H}+|C_2|^2\|Y^{[2]}\|^2_\mathbb{H}\ge c_{11}\left(\frac{|C_1|^2}{|{\rm Im}\lambda|+1}+\frac{|C_2|^2e^{2\pi|{\rm Im}\lambda|}}{|{\rm Im}\lambda|+1}\right)\eqno(12)
$$
$(c_{11}>0)$.

Obviously,
$$
\|C_1Y^{[1]}+C_2Y^{[2]}\|^2_\mathbb{H}=|C_1|^2\|Y^{[1]}\|^2_\mathbb{H}+|C_2|^2\|Y^{[2]}\|^2_\mathbb{H}+2{\rm Re}\langle C_1Y^{[1]},C_2Y^{[2]}\rangle\eqno(13)
$$

It follows from (2), (3) that

$$
\begin{array}{c}
|{\rm Re}\langle C_1Y^{[1]},C_2Y^{[2]}\rangle|\le|C_1C_2\int_0^\pi((1+b_{+,11}(x,\lambda))\overline{\elx b_{+,12}(x,\lambda)\elxx} \\\\+
b_{+,21}(x,\lambda)\elx\overline{(1+ b_{+,22}(x,\lambda))\elxx})dx|=
|C_1||C_2|\int_0^\pi o(1)e^{2i{\rm Re}\lambda x}dx=o(1)|C_1||C_2|.
\end{array}\eqno(14)
$$
From well-known inequality
$$
\varepsilon a^2+b^2/\varepsilon\ge2ab,
$$
where $a,b\ge0, \varepsilon>0$, taken for
$$
\varepsilon=\frac{|{\rm Im}\lambda|+1}{e^{2\pi|{\rm Im}\lambda|}}, \quad a=|C_1|, \quad b=|C_2|,
$$
 and (14) it follows that
$$
 2|C_1||C_2|\le \frac{|C_1|^2(|{\rm Im}\lambda|+1)}{e^{2\pi|{\rm Im}\lambda|}}+\frac{|C_2|^2e^{2\pi|{\rm Im}\lambda|}}{|{\rm Im}\lambda|+1}\le\frac{|C_1|^2}{|{\rm Im}\lambda|+1}+\frac{|C_2|^2e^{2\pi|{\rm Im}\lambda|}}{|{\rm Im}\lambda|+1}.\eqno(15)
$$
By virtue of (12), (14), and (15) for all sufficiently large $|{\rm Im}\lambda|$

$$
2|{\rm Re}\langle C_1Y^{[1]},C_2Y^{[2]}\rangle|\le(|C_1|^2\|Y^{[1]}\|^2_\mathbb{H}+|C_2|^2\|Y^{[2]}\|^2_\mathbb{H})/2,\eqno(16)
$$
hence, for all sufficiently large  $|{\rm Im}\lambda|$ it follows from (13), (16) that
$$
\|C_1Y^{[1]}+C_2Y^{[2]}\|^2_\mathbb{H}\ge(|C_1|^2\|Y^{[1]}\|^2_\mathbb{H}+|C_2|^2\|Y^{[2]}\|^2_\mathbb{H})/2.\eqno(17)
$$
Relations (12) and (17) imply
$$
\|C_1Y^{[1]}+C_2Y^{[2]}\|^2_\mathbb{H}\ge c_{12}\left(\frac{|C_1|^2}{|{\rm Im}\lambda|+1}+\frac{|C_2|^2e^{2\pi|{\rm Im}\lambda|}}{|{\rm Im}\lambda|+1}\right),\eqno(18)
$$

$$
\|C_1Y^{[1]}+C_2Y^{[2]}\|_\mathbb{H}\ge c_{12}(|C_1|\|Y^{[1]}\|_\mathbb{H}+|C_2|\|Y^{[2]}\|_\mathbb{H})\eqno(19)
$$
$(c_{12}>0)$.

Reasoning as above we find

$$
\frac{\tilde c_9e^{\pi|{\rm Im}\lambda|}}{\sqrt{|{\rm Im}\lambda|}+1}\le||Z^{[1]}||_\mathbb{H}\le\frac{\tilde c_{10}e^{\pi|{\rm Im}\lambda|}}{\sqrt{|{\rm Im}\lambda|}+1},\eqno(11^*)
$$

$$
\frac{\tilde c_9}{\sqrt{|{\rm Im}\lambda}|+1}\le||Z^{[2]}||_\mathbb{H}\le\frac{\tilde c_{10}}{\sqrt{|{\rm Im}\lambda|}+1},\eqno(10^*)
$$

$$
\|\tilde C_1Z^{[1]}+\tilde C_2Z^{[2]}\|^2_\mathbb{H}\ge \tilde c_{12}\left(\frac{|\tilde C_1|^2e^{2\pi|{\rm Im}\lambda|}}{|{\rm Im}\lambda|+1}+\frac{|\tilde C_2|^2}{|{\rm Im}\lambda|+1}\right),\eqno(18^*),
$$
and

$$
\|\tilde C_1Z^{[1]}+\tilde C_2Z^{[2]}\|_\mathbb{H}\ge  \tilde c_{12}(|\tilde C_1|\|Z^{[1]}\|_\mathbb{H}+|\tilde C_2|\|Z^{[2]}\|_\mathbb{H})\eqno(19^*)
$$

It follows from (18), (18*) that
$$
\begin{array}{c}
\|C_1Y^{[1]}+C_2Y^{[2]}\|^2_\mathbb{H}\|\tilde C_1Z^{[1]}+\tilde C_2Z^{[2]}\|^2_\mathbb{H}\\\\\ge c_{14}
\left(\frac{|C_1|^2}{|{\rm Im}\lambda|+1}+\frac{|C_2|^2e^{2\pi|{\rm Im}\lambda|}}{|{\rm Im}\lambda|+1}\right)\left(\frac{|\tilde C_1|^2e^{2\pi|{\rm Im}\lambda|}}{|{\rm Im}\lambda|+1}+\frac{|\tilde C_2|^2}{|{\rm Im}\lambda|+1}\right),
\end{array}
$$
therefore,
$$
\begin{array}{c}
\|C_1Y^{[1]}+C_2Y^{[2]}\|_\mathbb{H}\|\tilde C_1Z^{[1]}+\tilde C_2Z^{[2]}\|_\mathbb{H}\\\\\ge c_{15}
\left(\frac{|C_1|}{\sqrt{|{\rm Im}\lambda|}+1}+\frac{|C_2|e^{\pi|{\rm Im}\lambda|}}{\sqrt{|{\rm Im}\lambda|}+1}\right)\left(\frac{|\tilde C_1|e^{\pi|{\rm Im}\lambda|}}{\sqrt{|{\rm Im}\lambda|}+1}
+\frac{|\tilde C_2|}{\sqrt{|{\rm Im}\lambda|}+1}\right)\\\\
=c_{15}\left(
\frac{|C_1||\tilde C_1|e^{\pi|{\rm Im}\lambda|}}{{(\sqrt{|{\rm Im}\lambda|}+1})^2}+\frac{|C_1||\tilde C_2|}{{(\sqrt{|{\rm Im}\lambda|}+1})^2}+
\frac{|C_2||\tilde C_1|e^{2\pi|{\rm Im}\lambda|}}{{(\sqrt{|{\rm Im}\lambda|}+1})^2}+\frac{|C_2||\tilde C_2|e^{\pi|{\rm Im}\lambda|}}{{(\sqrt{|{\rm Im}\lambda|}+1})^2}\right).
\end{array}\eqno(20)
$$

By virtue of (19), (19*)
$$
\|C_1Y^{[1]}+C_2Y^{[2]}\|_\mathbb{H}\|\tilde C_1Z^{[1]}+\tilde C_2Z^{[2]}\|_\mathbb{H}\ge c_{16}\sum_{j,k=1}^2|C_j||\tilde C_k|\|Y^{[j]}\|_\mathbb{H}\|Z^{[k]}\|_\mathbb{H}.\eqno(21)
$$

Let us estimate the inner product $\langle\mathbf{y},\mathbf{z} \rangle$. Obviously,
$$
\begin{array}{c}
\langle C_1Y^{[1]}+C_2Y^{[2]},\tilde C_1Z^{[1]}+\tilde C_2Z^{[2]}\rangle=C_1\overline{\tilde C_1}\langle Y^{[1]},Z^{[1]}\rangle+C_1\overline{\tilde C_2}\langle Y^{[1]},Z^{[2]}\rangle\\\\+
C_2\overline{\tilde C_1}\langle Y^{[2]},Z^{[1]}\rangle+C_2\overline{\tilde C_2}\langle Y^{[2]},Z^{[2]}\rangle
\end{array}
$$
It follows from (2), (3), (2*), (3*) that
$$
\begin{array}{c}
\langle Y^{[1]},Z^{[1]}\rangle=\int_0^\pi((1+b_{+,11}(x,\lambda))\elx\overline{(1+\tilde b_{-,11}(x,\bar\lambda))\elxb}\\\\+b_{+,21}(x,\lambda)\elx\overline{\tilde b_{-,21}(x,\bar\lambda)\elxb})dx=\int_0^\pi(1+o(1))dx=O(1),
\end{array}\eqno(22)
$$

$$
\begin{array}{c}
\langle Y^{[1]},Z^{[2]}\rangle=\int_0^\pi((1+b_{+,11}(x,\lambda))\elx\overline{\tilde b_{-,12}(x,\bar\lambda)\elxxb}\\\\+b_{+,21}(x,\lambda)\elx\overline{(1+ \tilde b_{-,22}(x,\bar\lambda))\elxxb})dx=\int_0^\pi o(1)e^{2i\lambda x}dx=\frac{o(1)}{|{\rm Im}\lambda|+1},
\end{array}\eqno(23)
$$

$$
\begin{array}{c}
\langle Y^{[2]},Z^{[1]}\rangle=\int_0^\pi( b_{+,12}(x,\lambda)\elxx \overline{(1+\tilde b_{-,11}(x,\bar\lambda))\elxb}\\\\+(1+ b_{+,22}(x,\lambda))\elxx \overline{\tilde b_{-,21}(x,\bar\lambda)\elxb})dx=\int_0^\pi o(1)e^{-2i\lambda x}dx=\frac{o(1)e^{2\pi|{\rm Im}\lambda|}}{|{\rm Im}\lambda|+1},
\end{array}\eqno(24)
$$

$$
\begin{array}{c}
\langle Y^{[2]},Z^{[2]}\rangle=\int_0^\pi (b_{+,12}(x,\lambda)\elxx \overline{b_{-,12}(x,\bar\lambda)\elxxb }\\\\+(1+ b_{+,22}(x,\lambda))\elxx \overline{(1+ \tilde b_{-,22}(x,\bar\lambda))\elxxb })dx=\int_0^\pi (1+o(1))dx=O(1).
\end{array}\eqno(25)
$$

Combining (10), (11), (10*), (11*), (22-25), we obtain
$$
\langle Y^{[j]},Z^{[k]}\rangle=o(1)\|Y^{[j]}\|_\mathbb{H}\|Z^{[k]}\|_\mathbb{H}.
$$
$(j,k=1,2)$. The last relation implies
$$
\begin{array}{c}
|\langle\mathbf{y},\mathbf{z} \rangle|=
|\langle C_1Y^{[1]}+C_2Y^{[2]},\tilde C_1Z^{[1]}+\tilde C_2Z^{[2]}\rangle|\le \sum_{j,k=1}^2|C_j||\tilde C_k||\langle Y^{[j]},Z^{[k]}\rangle|\\\\=
 o(1)\sum_{1\le j,k\le1}|C_j||\tilde C_k|\|Y^{[j]}\|_\mathbb{H}\|Z^{[k]}\|_\mathbb{H}\\\\=o(1)(|C_1|\|Y^{[1]}\|_\mathbb{H}+|C_2|\|Y^{[2]}\|_\mathbb{H})(|\tilde C_1|\|Z^{[1]}\|_\mathbb{H}+|\tilde C_2|\|Z^{[2]}\|_\mathbb{H}).
\end{array}\eqno(26)
$$
 It follows from (19), (19*), (26) that
$$
|\langle\mathbf{y},\mathbf{z} \rangle|=o(1)\|\mathbf{y}\|_\mathbb{H}\|\mathbf{z}\|_\mathbb{H},
$$
hence,

$$
\|\mathbf{y}\|_\mathbb{H}\|\mathbf{z}\|_\mathbb{H}\to\infty
$$
as $|{\rm Im}\lambda|\to\infty$ if
$$
\langle\mathbf{y},\mathbf{z} \rangle=1.
$$
The case $\lambda\in\Omega_-$ can be considered by the similar way.

{\bf Theorem 1.} {\it Suppose the following conditions are valid:

1) Let $\{\mathbf{y}_n\}={\rm col}(y_{1,n}(x),y_{2,n}(x))$  be an arbitrary complete in $\mathbb{H}$ and minimal system consisting of eigen- and associated functions of operator  $\mathbb{L}$, and $\lambda_n$ be the corresponding eigenvalues $(n\in\mathbb{N})$;

2) Let  the biorthogonally adjoint system $\{\mathbf{z}_n\}={\rm col}(z_{1,n}(x),z_{2,n}(x))$ $(\langle\mathbf{y}_n,\mathbf{z}_k \rangle=\delta_{n,k})$ consist of eigen- and associated functions of operator  $\mathbb{L}^*$, and $\bar\lambda_n$ be the corresponding eigenvalues;

3) The systems  $\{\mathbf{y}_n\}$ and  $\{\mathbf{z}_n\}$ contain a finite number of associated functions;

4) $\overline{\lim}_{n\to\infty}|{\rm Im}\lambda_n|\to\infty$.

Then, the system $\{\mathbf{y}_n\}$ is not a basis in $\mathbb{H}$.
}

Proof.
Following [7], for arbitrary elements ${\bf y},{\bf z}$ from $\mathbb{H}$
 we introduce operator $\mathbf{y}\mathbf{z}$ in $\mathbb{H}$, acting according to the rule $\mathbf{f}\longrightarrow\langle\mathbf{f},\mathbf{z}\rangle\mathbf{y}$,
or in the scalar form

$$
\langle\mathbf{f},\mathbf{z}\rangle\mathbf{y}=
\begin{pmatrix}
(\int_0^\pi(f_1(t)\overline{z_1(t)}+f_2(t)\overline{z_2(t)})dt)y_1(x)\\
(\int_0^\pi(f_1(t)\overline{z_1(t)}+f_2(t)\overline{z_2(t)})dt)y_2(x)
\end{pmatrix}.\eqno(27)
$$
 It is easy to see that the operator $\mathbf{y}\mathbf{z}$ is determined by the matrix ${\bf y}\overline{{\bf z}}^T$
with entries
$y_j(x)\overline{z_k(t)}$ $(1\le j,k\le2)$.
Let us define the norm of the indicated operator by the equality
 $$
 \|\mathbf{y}\mathbf{z}\|=\sup_{\|\mathbf{f}\|_\mathbb{H}\le1}\|\langle\mathbf{f},\mathbf{z}\rangle\mathbf{y}\|_\mathbb{H}.
$$
Let us estimate $\|\mathbf{y}\mathbf{z}\|$. Let $\|z_2\|\ne0$. Set $\mathbf{\tilde f}={\rm col}(\tilde f_1,\tilde f_2))$, where $\tilde f_1=0, \tilde f_2=z_2/\|z_2\|$, then, according to (27)
$$
\|\langle\mathbf{\tilde f},\mathbf{z}\rangle\mathbf{y}\|=\|z_2\|\|\mathbf{y}\|_\mathbb{H}\ge\|z_2\|(\|y_1\|+\|y_2\|)/\sqrt{2}.\eqno(28)
$$
If  $\|z_2\|=0$,  the last inequality is trivial.

Let $\|z_1\|\ne0$. Set $\mathbf{\hat f}={\rm col}(\hat f_1,\hat f_2))$, where $\hat f_2=0, \hat f_1=z_1/\|z_1\|$, then, according to (27)
$$
\|\langle\mathbf{\hat f},\mathbf{z}\rangle\mathbf{y}\|=\|z_1\|\|\mathbf{y}\|_\mathbb{H}\ge\|z_1\|(\|y_1\|+\|y_2\|)/\sqrt{2}.\eqno(29)
$$
If  $\|z_1\|=0$, the last inequality is trivial. It follows from (28), (29) that
$$
\|\mathbf{y}\mathbf{z}\|\ge(\|z_1\|+\|z_2\|)(\|y_1\|+\|y_2\|)/(2\sqrt{2})\ge\|\mathbf{y}\|_\mathbb{H}\|\mathbf{z}\|_\mathbb{H}/(2\sqrt{2}).\eqno(30)
$$

For all $n>n_0$ $\mathbf{y}_n$ and $\mathbf{z}_n$ are eigenfunctions.
Denote by
$ \mathfrak{P}_n$ the operator determined by the  matrix
$$
\begin{pmatrix}
y_{1,n}(x)\overline{ z_{1,n}(t)}&y_{1,n}(x)\overline{ z_{2,n}(t)}\\
y_{2,n}(x)\overline{ z_{1,n}(t)}&y_{2,n}(x)\overline{ z_{1,n}(t)}
\end{pmatrix}
$$
and acting according to  formula (27). It follows from (30) and Lemma 1 that
$$
\overline{\lim}_{n\to\infty}\|\mathfrak{P}_n\|=\infty.
$$
This, together with a theorem of resonance type  [8, Ch. 2, p. 1] implies existence of a function $\mathbf{f}\in\mathbb{H}$ such that
$$
\overline\lim_{|n|\to\infty}|\langle\mathbf{f},\mathbf{z}_n\rangle|\|\mathbf{y}_n\|_\mathbb{H}=\infty,
$$
hence, the root function system $\{\mathbf{y}_n\}$  is not a basis in $\mathbb{H}$.

This completes the proof.
\medskip
\medskip

{\bf Remark 1.} A lot of boundary value problems for operator (1)
 with two-point boundary conditions
$$
\begin{array}{c}
U_1(\mathbf{y})= a_{11}y_1(0)+a_{12}y_2(0)+a_{13}y_1(\pi)+ a_{14}y_2(\pi)=0,\\ U_2(\mathbf{y})= a_{21}y_1(0)+a_{22}y_2(0)+a_{23}y_1(\pi)+ a_{24}y_2(\pi)=0
\end{array}\eqno(31)
$$
satisfies all the conditions of Theorem 1, for instance, a problem for the Dirac operator with boundary conditions of type (31) studied by A.P. Kosarev and A.A. Shkalikov in [2].

\medskip
\medskip
\medskip
\medskip

\centerline {References}
\medskip
\medskip

[1] A. M. Savchuk and I. V. Sadovnichaya.
Spectral Analysis of 1D Dirac System with Summable
Potential and Sturm-Liouville Operator
with Distribution Coefficients. Differential Equations, 2024, Vol. 60, Suppl. 2, pp. 145-348.

[2] A.P. Kosarev,  A.A. Shkalikov. Spectral Asymptotics of Solutions of a $2\times 2$ System of First-Order Ordinary Differential Equations. Math Notes, 2021, Vol. 110, 967-971 .

[3] A.P. Kosarev,  A.A. Shkalikov. Asymptotics in the Spectral Parameter for Solutions of $2\times 2$ Systems of Ordinary Differential Equations.  Math Notes, 2023, Vol. 114, No. 4, 472-488 .

[4] A.P. Kosarev,  A.A. Shkalikov. Asymptotic Expansions of Solutions to $n\times n$ Systems of Ordinary Differential Equations with a Large Parameter. Math Notes, 2024, Vol.116, No. 6, 1312-1325 .

[5] A.P. Kosarev,  A.A. Shkalikov. Asymptotic Representations of Solutions of $n\times n$ Systems  of Ordinary Differential Equations with a Large Parameter. Math Notes, 2024,  Vol. 116, No. 2, 283-302.

[6] M.M. Malamud  and  L.L. Oridoroga. On the completeness of root subspaces of boundary value problems
for first order systems of ordinary differential equations, J. Funct. Anal., 2012, Vol. 263, pp. 1939-1980.

[7] M. V. Keldysh. On the completeness
of the eigenfunctions of some classes of non-selfadjoint
linear operators,
Russian Mathematical Surveys, 1971, Vol. 26, Issue 4, pp.
15–44.

[8] K. Yosida. Functional Analysis. — Springer, 1980.
\medskip
{\normalsize

\medskip
\medskip
\medskip

email: alexmakin@yandex.ru

\end{document}